\documentclass[12pt]{amsart}
\usepackage{graphicx,amssymb,amsmath}
\newtheorem{thm}{Theorem}[section]

\newtheorem{lem}[thm]{Lemma}
\newtheorem{prop}[thm]{Proposition}
\theoremstyle{definition}

\theoremstyle{remark}

\numberwithin{equation}{section}

\newcommand{\real}{\mathbb R}

\begin{document}
\title{Enhancement of combustion by drift in a coupled reaction-diffusion model}
\author{Lam Raga A. Markely, David Andrzejewski, Erick Butzlaff and Alexander Kiselev }
\thanks{Department of
Mathematics, University of Wisconsin, Madison, WI 53706; e-mail:
kiselev@math.wisc.edu. This work has been a part of Collaborative
Undergraduate Research Lab at the University of Wisconsin. Support
of the NSF VIGRE grant is acknowledged.}

\begin{abstract}
We study analytically and numerically a model describing front
propagation of a KPP reaction in a fluid flow. The model consists
of coupled one-dimensional reaction-diffusion equations with
different drift coefficients. The main rigorous results give lower
bounds for the speed of propagation that are linear in the drift
coefficient, which agrees very well with the numerical
observations. In addition, we find the optimal constant in a
functional inequality of independent interest used in the proof.
\end{abstract}
\maketitle

\section{Introduction}

Many reaction processes in nature and engineering take place in a
moving fluid.  The creation of ozone in the atmosphere, the
nuclear explosion of a supernova, pattern formation in
morphogenesis, wild fires, and the gasoline transformation in an
internal combustion engine are just a few examples. The effect of
fluid flow on a reaction process can be profound, especially if
the flow is strong. The problem has been studied for many years by
engineers, physicists, and mathematicians alike.
One of the most extensively studied mathematical models of the
premixed reaction process is the reaction-diffusion equation and
systems. The advection-reaction-diffusion equation is given by
\begin{equation}\label{rd1}
T_t + Au\cdot\nabla T = \kappa\Delta T +  M f(T).
\end{equation}
Here $u(x)$ is the fluid flow, which we will assume is
time-independent and passive. The coefficient $A$ is the flow
strength parameter, $\kappa$ is the diffusivity, $M$ is the
reaction strength, and $f(T)$ is the reaction term. The function
$T(x,t)$ is normalized so that $0 \leq T(x,t) \leq 1,$ and can
stand for normalized temperature, the mass fraction of a reactant,
or a share of a population with a certain feature, depending on the problem.
Classical works of Fisher and Kolmogorov-Petrovski-Piskunov (KPP)
\cite{Fisher,KPP} first considered equation \eqref{rd1} in the dimension
$d=1$ and with $u=0$. These works modeled the propagation of an advantageous
gene in the population, and established existence and, in a certain sense,
stability of the traveling fronts. The extensions have occupied some of
the best efforts of mathematicians since. Equation \eqref{rd1} in several
spacial dimensions in the absence of advection ($u=0$) is by now fairly well
understood: the existence of traveling waves, stability, and asymptotic propagation properties have been
extensively studied (see, e.g. \cite{AW,freidlin-4,HN,VVV} where
many more references can be found). The effects of advection,
however, are very important in many situations in combustion and
more general chemical reactions \cite{P1,Wil,ZBLM} as well as in
biology and ecology \cite{Levin,Murray}.

One important phenomena that is observed in many situations is the
enhancement of the reaction rate by fluid motion. The physical
reason for this observed speed-up is believed to be that fluid
advection tends to increase the area available for reaction.
Recently, there has been significant progress in the mathematical
understanding of the reaction rate enhancement for several different
classes of flows
\cite{ACVV1,ACVV2,ABP,CKOR,CKR,HPS,KB,KBM,KS,KRS,KR,MS,NX,VCKRR}.
The papers \cite{Berrev,Xin2} provide excellent reviews of some of
these results and further references.

In this paper, our goal is to study reaction enhancement within a
framework of a model given by a coupled system of one-dimensional
advection-reaction-diffusion equations:
\begin{equation}\label{mainmodel}
(T_{j})_t + A_j (T_j)_x - (T_j)_{xx} =  T_j(1-T_j) +
\alpha(T_{j-1}-T_j) + \alpha(T_{j+1}-T_j),
\end{equation}
$j =1, \dots, N.$ We assume that the model is periodic in $j;$
that is, in the equation for $T_1(x),$ we replace $T_0(x)$ with
$T_N(x),$ and in the equation for $T_N(x),$ we put $T_1(x)$
instead of $T_{N+1}(x).$ The initial data always satisfies $0
\leq T_j(x,0) \leq 1,$ and it is a simple consequence of the
maximum principle that $T_j(x,t)$ remains within these bounds for
all times (see Section~\ref{aux}). Our reaction
term is the classical KPP \cite{KPP}, and
 for simplicity we set $\kappa =
M=1$ (these parameters are normalized by a simple rescaling of
time and space). We assume that $\sum_{j=1}^N A_j =0,$ which
corresponds to the mean zero flow. Any nonzero mean is taken into
account by switching to a moving system of coordinates; we are
interested in a non trivial effect of front stretching rather than
simple transfer with a constant speed. One can regard
\eqref{mainmodel} as a model for reaction in a layered fluid,
where different layers move with different speeds. The role of the
parameter $\alpha>0$ is to provide diffusive coupling between
different layers, and its meaning will be further discussed in
Section~\ref{main}. Our main interest is in the case of the large
coefficients $A_j$, where the influence of the drift is most
pronounced. We consider the front-like initial data where
$T_j(x,0) \rightarrow 0$ as $x \rightarrow \infty$ and $T_j(x,0)
\rightarrow 1$ as $x \rightarrow -\infty$ exponentially fast (we
will state the exact conditions in Section~\ref{aux}). We adopt
the following natural quantity (called the bulk burning rate in
\cite{CKOR}) as the main measure of the reaction rate:
\begin{equation}\label{bbrdef}
V(t) = \frac{1}{N}\sum\limits_{j=1}^N
\int\limits_{-\infty}^\infty (T_j)_t(x) \,dx.
\end{equation} One of
the results that we prove is as follows.

\begin{thm}\label{main1}
Assume that the initial data $T_j(x,0)$ are front like
\eqref{id1}.
Then there
exists a universal constant $C$ such that for any $\tau>0,$ we
have
\begin{equation}\label{lb1}
\frac{1}{\tau}\int\limits^{\tau}_{0}V(t) dt \geq
\frac{C}{N}\sum\limits_{j=1}^N|A_{j-1}-A_j|
\left[\frac{1}{\alpha\tau^2}+ \frac{1}{\alpha\tau}+\frac{1}{\tau}+
\alpha+1\right]^{-1}.
\end{equation}
\end{thm}

There are several noteworthy properties the lower bound
\eqref{lb1} possesses. Firstly, the bound is linear in the flow
strength. Secondly, the estimate clearly shows the need to wait
for a certain time $\tau$ before the front propagation (and thus
the reaction rate) stabilizes and the long-time lower bound is
achieved. If the coupling constant $\alpha$ is large, we have the
lower bound $CA \alpha^{-1}$ starting from times $\tau \approx 1.$
If $\alpha$ is small, the lower bound is $CA,$ starting from
times $\tau \approx \alpha^{-1}.$ We interpret these results
further in Sections~\ref{main} and \ref{num}. The numerical
experiments discussed in Section~\ref{num} show very good
qualitative agreement with \eqref{lb1}.

We remark that the result of Theorem~\ref{main1} is reminiscent of
the bounds proved in \cite{CKOR} for the two dimensional equation
\eqref{rd1} with a shear flow $u = (u(y),0).$ This is not
surprising since our system \eqref{mainmodel} can be regarded as a
discrete model of a shear flow. However, the proofs in the case of
model \eqref{mainmodel} are conceptually much more transparent,
the estimates are more precise, and the dependence on the
averaging time is more explicit.

The paper is organized as follows. In Section~\ref{aux} we collect
some background estimates which are needed in the proof of main
results. In particular, we prove a sharp version of a functional
inequality which was first discovered in \cite{CKOR}. The
inequality is somewhat reminiscent of the indeterminacy principle
of quantum mechanics. Although finding a sharp constant is not
particularly important for the key results of this paper, we feel
that the result is elegant (even if elementary) and may be of
independent interest. In Section~\ref{two} we discuss a simpler
two-layer case. In Section~\ref{main} we prove
Theorem~\ref{main1}, and some related results for a slightly more
general model. In Section~\ref{num} we describe the numerical
simulations.

\section{The Auxiliary Tools and The Optimal Constant}\label{aux}

We make our model slightly more general by associating a width
$h_j$ with each layer. By dimensionality of the discrete gradient
terms providing coupling of the neighboring layers in
\eqref{mainmodel}, the parameter $\alpha$ is made to depend on $j,$ with $\alpha_j =\kappa /h_j:$
\begin{equation}\label{mainmodel1}
(T_{j})_t + A_j (T_j)_x - (T_j)_{xx} =  T_j(1-T_j) +
\alpha_j(T_{j-1}-T_j) + \alpha_j(T_{j+1}-T_j).
\end{equation}
The parameter $\kappa$ now plays the role of diffusivity between the
layers. As opposed to the diffusivity in $x$ direction, this
diffusion coefficient cannot be set to one by a simple rescaling.
Our mean zero flow condition now reads $\sum_{j=1}^N h_jA_j =0.$ Set
$H=\sum_{j=1}^N h_j.$

We first prove a lemma on the preservation of spacial decay of the
solutions of \eqref{mainmodel1} that will allow us to manipulate
the equations, in particular integrate by parts.

\begin{lem}\label{decay}
Assume that $T_j(x,t),$ $j=1,\dots,N$ satisfy \eqref{mainmodel1}
and that the initial data $T_j(x,0)$ satisfy
\begin{equation}\label{id1}
T_j(x,0) < C_0 e^{-\lambda x}, \,\,\,\,1-T_j(x,0) < C_0 e^{\lambda
x},\,\,\,\,|(T_j(x,0))_{x}| < C_0 e^{-\lambda |x|}
\end{equation}
for all $j$ and some $\lambda>0.$ Assume that $c_1$ is such that
\begin{equation}\label{upperspeed}
c_1 \geq {\rm max}_j |A_j| +\lambda+ \lambda^{-1}(1+ 2{\rm max}_j
\alpha_j).
\end{equation}
Then for all $x,$ $t > 0$ we have
\begin{equation}\label{expdecay}
T_j(x,t) \leq C_0 e^{-\lambda(x-c_1t)}, \,\,\,1-T_j(x,t) < C_0
e^{\lambda (x+c_1t)},\,\,\,|(T_j(x,t))_{x}| < C_0 e^{\lambda
(|x|+c_1t)}.
\end{equation}
\end{lem}
\begin{proof}
Let us start by establishing the first bound in \eqref{expdecay};
the second bound can be proved by an identical argument. Set
$\phi(x,t) = e^{-\lambda (x-c_1t)},$ and $G_j(x,t) = C_0 \phi(x,t)
- T_j(x,t).$ A direct computation using \eqref{upperspeed} shows
that
\begin{equation}\label{phi} \phi_t + max_{j}|A_j| \phi_x -\phi_{xx}
-(1+2{\rm max}_j \alpha_j)\phi \geq 0 \end{equation} for any $j.$
Note that $G_j(x,0)>0$ in view of \eqref{id1}. On the contrary, assume
that the first bound in \eqref{expdecay} is not satisfied
and that $t_0>0$ is the first time such that there exists $j$ and
$x_0$ so that $G_j(x_0,t_0)=0.$ Combining \eqref{phi} with
\eqref{mainmodel1}, we find that
\begin{equation}\label{subsol}
(G_j)_t +A_j (G_j)_x -(G_j)_{xx} \geq 0
\end{equation} for all $x,$ $t \leq t_0.$ Now, in \eqref{subsol},
$G_j(x,0)>0$ for all $x$ and $G_j(x_0,t_0)=0$ contradict the
well-known maximum principle for parabolic PDE (see e.g.
\cite{Evans}). The proof of the second inequality in
\eqref{expdecay} is similar and is omitted. To prove the
bound for $(T_j)_x,$ differentiate \eqref{mainmodel1}; denoting
$H_j = (T_j)_x,$ we have
\[ (H_j)_t + A_j (H_j)_x - (H_j)_{xx} =H_j (1-2T_j) +
\alpha_j(H_{j-1}-H_j) + \alpha_{j+1}(H_{j+1}-H_j). \] Given that
$0 \leq T_j(x)\leq 1$ by a simple application of maximum
principle, the bounds for $(T_j)_x$ are now obtained in the same
manner as above.
\end{proof}

The bulk burning rate for the model \eqref{mainmodel1} is now
defined accordingly via
\begin{equation}\label{bbr1}
V(t) = \frac{1}{H}\sum\limits_{j=1}^N
h_j\int\limits_{-\infty}^\infty (T_j)_t(x) \,dx.
\end{equation}
Assuming the initial data satisfies \eqref{id1}, the indefinite
integrals are well-defined (integrating over $[-B,B]$ and taking
$B$ to infinity). As a consequence of \eqref{mainmodel1},
Lemma~\ref{decay}, and the condition $\sum_{j=1}^N h_j A_j =0,$ we also
have
\begin{equation}\label{bbr2}
V(t) = \frac1H\sum\limits_{j=1}^N h_j\int\limits_{-\infty}^\infty
T_j(x)(1-T_j(x)) \,dx.
\end{equation}
We continue by computing:
\begin{eqnarray}
V' (t) = \partial_t \frac1H\sum\limits_{j=1}^N
h_j\int\limits_{-\infty}^\infty T_j(x)(1-T_j(x)) \,dx =
\frac1H\sum\limits_{j=1}^N h_j\int\limits_{-\infty}^\infty
(T_j)_t(x)(1-2T_j(x)) \,dx \nonumber \\
= \frac1H\sum\limits_{j=1}^N h_j\int\limits_{-\infty}^\infty
\left(2((T_j)_{x})^2 + T_j(1-T_j)(1-2T_j(x)) +
\frac{2\kappa}{h_j}(T_j(x) - T_{j-1}(x))^2\right) \,dx \nonumber \\
\geq \frac1H\sum\limits_{j=1}^N h_j\int\limits_{-\infty}^\infty
\left(2((T_j)_{x})^2  + \frac{2\kappa}{h_j}(T_j(x) -
T_{j-1}(x))^2\right) \,dx - V(t) \label{keycontrol}.
\end{eqnarray}
In the second step above, we substitute $(T_j)_t(x)$ from
\eqref{mainmodel1}, integrate by parts, and rearrange the terms.
These manipulations are justified by Lemma~\ref{decay}. In the last
step we use the expression \eqref{bbr2} and the fact that $0 \leq
T_j(x) \leq 1.$

As a warm up before our main results, we prove the following
estimate.
\begin{thm}\label{general}
For any choice of $A_j$ such that $\sum_j h_j A_j =0,$ we have
\begin{equation}\label{genbound}
V(t)^2 \geq \frac{\pi^2}{16} + e^{-2t}\left( V(0)^2 -
\frac{\pi^2}{16} \right).
\end{equation}
\end{thm}
A key element in the proof is a general functional inequality.
\begin{prop}\label{fineq}
For any function $T(x) \in C^1(\real),$ $0 \leq T(x) \leq 1,$ such
that $T(x) \to 1$ as $x\to -\infty$ and $T(x) \to 0$ as $x\to
\infty,$ we have
\begin{equation}\label{funin}
\int\limits_{-\infty}^\infty
\left(T_x\right)^2dx\int\limits_{-\infty}^\infty T(1-T)dx \geq
\left(\frac{\pi}{8}\right)^2.
\end{equation}
The constant $\left(\frac{\pi}{8}\right)^2$ is sharp.
\end{prop}
\begin{proof}
We note that the inequality \eqref{funin} has been proved in
\cite{CKOR} with a weaker universal constant on the right hand
side. Using the Cauchy-Schwartz inequality we find that
\begin{eqnarray*}
\int\limits_{-\infty}^\infty\left(T_x\right)^2dx\int\limits_{-\infty}^\infty
T(1-T)dx &\geq& \left(\,\int\limits_{-\infty}^\infty
T_x\left(T(1-T)\right)^{1/2}dx\right)^2.
\end{eqnarray*}
Changing the variable and taking into account the asymptotic
behavior of $T(x)$, we have
\[\int\limits^{\infty}_{-\infty}T_x(T(x)(1-T(x)))^{1/2}dx =
\int\limits^{1}_{0}(T(1-T))^{1/2}dT.\] The latter integral is
computed explicitly by setting $T = \cos \theta$, and evaluates to
$\frac{\pi}{8}$. This proves \eqref{funin}. To show that this
constant is sharp, note that if $T_x = (T(1-T))^{1/2}$, then
\[\int\limits^{\infty}_{-\infty}\left(T_x\right)^2dx\int\limits^{\infty}_{-\infty}T(1-T)dx
= \left(\int\limits^{\infty}_{-\infty}T_x(T(1-T))^{1/2}dx\right)^2
= \left(\frac{\pi}{8}\right)^2. \qquad  \] We can solve for $T(x)$
explicitly and use the solution to construct an explicit function
satisfying the necessary asymptotic behavior for which the
equality in \eqref{funin} is achieved. In particular,

\[ T(x) = \left\{ \begin{array}{ll}
           1 &{\rm if} \,\,\, x \leq -\frac{\pi}{2} \\
\frac{1-\sin x}{2}
         & {\rm if} \,\,\, -\frac{\pi}{2}\leq x \leq \frac{\pi}{2} \\
           0         &{\rm if} \,\,\, \frac{\pi}{2} \leq
           x. \end{array} \right. \]
is one such function.
\end{proof}

\begin{proof}[Proof of Theorem~\ref{general}]
From \eqref{keycontrol}, we see that
\[ V' (t)+V(t) \geq \frac2H\sum\limits_{j=1}^N h_j\int\limits_{-\infty}^\infty
((T_j)_{x})^2dx. \] From \eqref{id1} and Proposition~\ref{fineq}, we
derive
\begin{eqnarray*}
 V' (t)+V(t) \geq  \frac{\pi^2}{32H}
\sum\limits_{j=1}^N \frac{h_j}{\int\limits_{-\infty}^\infty
T_j(x)(1-T_j(x))\,dx} \\ \geq \frac{\pi^2H}{32} \left(
\sum\limits_{j=1}^N h_j\int\limits_{-\infty}^\infty
T_j(x)(1-T_j(x))\,dx \right)^{-1}  = \frac{\pi^2}{32V(t)}.
\end{eqnarray*}
On the second step we use Cauchy-Schwartz inequality.
Therefore, we have
\[ \frac{d}{dt} (V(t)^2) +2V(t)^2 \geq \frac{\pi^2}{16}. \]
The inequality \eqref{genbound} follows from the standard
application of the Gronwall lemma.
\end{proof}
Theorem~\ref{general} gives a universal lower bound on the
reaction rate, independent of the flow. One consequence is that in
the framework of our model, a mean zero flow cannot quench the
reaction. However, we are interested in the bounds that show
enhancement of the reaction. We will prove such bounds in the next
two sections.

\section{The Two-Layer Model}\label{two}

For the sake of clarity, we start with a simplified version of our
general model \eqref{mainmodel1} where only two layers of unit
width are involved. Namely, we look at a system
\begin{eqnarray}\label{twolayer}
(T_1)_t + A(T_1)_t - (T_1)_{xx} = T_1(1-T_1) + \kappa (T_2 - T_1)\\
(T_2)_t - A(T_2)_t - (T_2)_{xx} = T_2(1-T_2) + \kappa (T_1 - T_2)
\nonumber
\end{eqnarray}
with the initial data $T_{1,2}(x,0)$ satisfying \eqref{id1}. We
set $h_1 =h_2 =1$ in the definition \eqref{bbr1} of the bulk
burning rate. Throughout the rest of the paper we denote by $C$
different universal constants which enter the estimates.

\begin{thm}\label{2layer}
There exists a universal positive constant $C$ such that for any
$\tau>0,$ we have
\begin{equation}\label{lb2l}
\frac{1}{\tau}\int\limits^{\tau}_{0}V(t) dt \geq CA
\left[\frac{1}{\kappa\tau^2}+ \frac{1}{\kappa\tau}+\frac{1}{\tau}+
\kappa+1\right]^{-1}.
\end{equation}
\end{thm}
\noindent \it Remark. \rm Note that Theorem~\ref{2layer} is a
particular case of Theorem~\ref{main1} for $N=2.$ \\

\begin{proof}
We first integrate \eqref{twolayer} over all real axis. Given the
asymptotic behavior of $T_1$ and $T_2,$ \eqref{id1} and
Lemma~\ref{decay}, we obtain
\begin{eqnarray}\label{t1}
& \int_\mathbb{R} (T_1)_t dx -A = \int_\mathbb{R}
T_1(1-T_1)dx + \kappa \int_\mathbb{R}(T_2 - T_1)dx\\
\label{t2} & \int_\mathbb{R}(T_2)_t dx +A = \int_\mathbb{R}
T_2(1-T_2)dx + \kappa \int_\mathbb{R}(T_1 - T_2)dx.
\end{eqnarray}
\noindent Subtracting \eqref{t2} from \eqref{t1} and integrating
the result over $t \in [\tau_1,\tau_2]$, we find that
\begin{equation}
\begin{split}\label{2l3}
&\int_\mathbb{R} (T_1(\tau_2)-T_2(\tau_2))dx -
\int_\mathbb{R}(T_1(\tau_1)-T_2(\tau_1))dx \\+
&\int\limits^{\tau_2}_{\tau_1} \int_\mathbb{R}
(T_2(1-T_2)-T_1(1-T_1))dxdt \\
&+ 2\kappa \int\limits^{\tau_2}_{\tau_1}
\int_\mathbb{R}(T_1-T_2)dx dt = 2 A (\tau_2-\tau_1).
\end{split}
\end{equation}
We now claim the following Lemma.
\begin{lem}\label{keylem}
For any $t>0,$ we have
\begin{eqnarray}\label{keyest0}
\int_\mathbb{R} | T_1(x,t)-T_2(x,t) | dx  \leq  3
\left[\int_\mathbb{R}(T_1(x,t)-T_2(x,t))^2 dx + \right. \\ \left.
\int_\mathbb{R}
\left(T_1(x,t)(1-T_1(x,t))+T_2(x,t)(1-T_2(x,t))\right)dx \right]
\leq \frac{3}{\kappa} \left( V(t)+V^\prime(t)\right) + 3 V(t).
\nonumber
\end{eqnarray}
\end{lem}
\begin{proof}
\noindent We begin by proving the first inequality, which is valid
for any functions $T_1(x),$ $T_2(x) \in L^1(\mathbb{R}) \cap
L^2(\mathbb{R}):$
\begin{equation}\label{keyest}
\int_\mathbb{R} | T_1-T_2 | dx \leq 3
\left[\int_\mathbb{R}(T_1-T_2)^2 dx + \int_\mathbb{R}
(T_1(1-T_1)+T_2(1-T_2))dx\right].
\end{equation}
Consider two different cases. Denote $S$ the set of all $x$ such
that $|T_1(x) - T_2(x)| \geq 1/3.$ Then $|T_1(x)-T_2(x)| \leq
3|T_1(x)-T_2(x)|^2$ and \eqref{keyest} holds with $C=3$ if we
restrict the integration to set $S.$ Next, assume that $|T_1(x)
-T_2(x)| \leq 1/3.$ Note that if at least one of $T_1(x),$ $T_2(x)$
belongs to the interval $[1/3,2/3],$ then
$T_1(x)(1-T_1(x))+T_2(x)(1-T_2(x)) \geq 2/9.$ Consider the
alternative where both $T_1(x)$ and $T_2(x)$ lie in either $[0,1/3)$
or $(2/3,1].$ Note that the absolute value of the derivative of
function $x(1-x)$ satisfies $|(1-2x)| \geq 1/3$ in these intervals.
Then by the mean value theorem,
\[ T_1(x)(1-T_1(x))+T_2(x)(1-T_2(x)) \geq |T_1(x)(1-T_1(x))-
T_2(x)(1-T_2(x))| \geq \frac13
|T_1(x) -T_2(x)|. \] Therefore, if $x \in \mathbb{R} \setminus S,$
then
\[ |T_1(x)-T_2(x)| \leq 3 \left(
T_1(x)(1-T_1(x))+T_2(x)(1-T_2(x))\right). \] Combining the two
cases, we obtain \eqref{keyest}.

Now the second inequality in \eqref{keyest0} follows from the
definition \eqref{bbr1} of bulk burning rate, $V(t),$ and the
estimate \eqref{keycontrol} (recall that we set $h_{1,2}=1$).
\end{proof}

Given Lemma~\ref{keylem}, \eqref{2l3} implies that
\begin{equation}\label{bbrred}
\frac{3}{\kappa} [V^\prime(\tau_2)+ V^\prime(\tau_1)]
+(\frac{3}{\kappa}+9)[V(\tau_2)+V(\tau_1)] + (6\kappa + 7)
\int\limits^{\tau 2}_{\tau 1} V(t)dt \geq 2A (\tau_2-\tau_1).
\end{equation}

Let us now
apply the following averaging over $t\in[0,\tau]$ to the difference
of equations \eqref{t1} and \eqref{t2}:

\begin{equation}\label{aver1}
\frac{1}{\tau^3} \int\limits^{\frac{\tau}{4}}_{0}
\int\limits^{\frac{\tau}{4}+s}_{\frac{\tau}{4}-s}
\int\limits^{\frac{\tau}{2}+k}_{\frac{\tau}{2}-k} dt dk ds =
\frac{1}{\tau^3} \int\limits^{\tau}_{0}
G(\frac{\tau}{2},t-\frac{\tau}{2})dt
\end{equation}
\noindent where
\begin{equation}\label{averker}
G(h,\xi) = \left\{ \begin{array}{ll} \frac{1}{2}(h-|\xi|)^2 -
(\frac{h}{2}-|\xi|)^2 & {\rm if} \,\,
 0 \leq |\xi|\leq \frac{h}{2}\\
 \frac{1}{2}(h-|\xi|)^2 & {\rm if} \,\, \frac{h}{2}\leq
|\xi| \leq h.
\end{array} \right.
\end{equation}

\noindent 
After the first integration, we obtain \eqref{2l3} and thus
\eqref{bbrred} with $\tau_1 = \tau-k,$ $\tau_2=\tau+k.$
\noindent Integrating twice more
leads to the following estimate:
\begin{eqnarray}\nonumber
\frac{3}{\kappa}\int\limits^{\tau}_{0}V(t)dt+\left(\frac{3}{\kappa}+9\right)
\int\limits^{\tau}_{0}V(t)h(t,\tau) dt  + \\  \label{2lf} (6\kappa
+7) \int\limits^{\tau}_{0} V(t)G(\frac{\tau}{2},t-\frac{\tau}{2})
dt \geq \frac{A\tau^3}{16}
\end{eqnarray}

\noindent where
\[ h(t,\tau) = \left\{ \begin{array}{ll} t
         & {\rm if} \,\,\, 0\leq t \leq \frac{\tau}{4} \\
           \frac{\tau}{4} &{\rm if} \,\,\, \frac{\tau}{4}\leq t \leq \frac{3\tau}{4} \\
           \tau-t         &{\rm if} \,\,\, \frac{3\tau}{4}\leq t \leq
           \tau. \end{array} \right. \]

\noindent Observe that $G(\frac{\tau}{2},t-\frac{\tau}{2}) \leq
\frac{\tau^2}{8}$ and that $h(t,\tau) \leq\frac{\tau}{4}$. Hence
\eqref{2lf} can be expressed as follows:
\[ \frac{1}{\tau}\int\limits^{\tau}_{0}V(t) dt \geq C
A\left[\frac{1}{\kappa\tau^2}+ \frac{1}{\kappa \tau}+
\frac{1}{\tau} + \kappa +1 \right]^{-1}. \]
This coincides with \eqref{lb2l}, thus proving the theorem. It is also
easy to estimate that the constant $C$ can be taken equal to
$1/48.$
\end{proof}

\section{The Multi-Layer Model}\label{main}

Here we prove a more general version of Theorem~\ref{main1}. Let
us introduce the notation
\[ R_1 =  1+ {\rm max}_j h_j^{-1} , \,\,\,R_2=  1+ {\rm
max}_j (h_j h_{j-1}^{-1}). \]

\begin{thm}\label{main2}
Assume that the initial data $T_j(x,0)$ for \eqref{mainmodel1} are
front like \eqref{id1}.
Then there exists a universal constant $C$ such that for any
$\tau>0,$ we have
\begin{equation}\label{lb2}
\frac{1}{\tau}\int\limits^{\tau}_{0}V(t) dt \geq
\frac{C}{HR_2}\sum\limits_{j=1}^N|A_{j-1}-A_j|h_j
\left[\frac{1}{\kappa\tau^2}+
\frac{1}{\kappa\tau}+\frac{R_1}{\tau}+ \kappa R_1+1\right]^{-1}.
\end{equation}
\end{thm}
\noindent \it Remark. \rm Theorem~\ref{main1} is a particular case
of Theorem~\ref{main2} when $h_j \equiv 1$ (with $\alpha$
substituted by $\kappa$). \\
\begin{proof}
Consider two neighboring layers, $j$ and $j-1.$ Assume without
loss of generality that $A_j \geq A_{j-1}.$ Subtracting the equation
for $T_j$ from the equation for $T_{j-1}$ and integrating, we find:
\begin{eqnarray}\label{diff1}
\int_\mathbb{R}((T_{j-1})_t-(T_j)_t)\, dx + A_j -A_{j-1} =
\int_\mathbb{R} \left( T_{j-1}(1-T_{j-1}) - T_j(1-T_j) \right)\,
dx + \\
\left(\frac{\kappa}{h_j}+\frac{\kappa}{h_{j-1}}\right)
\int_\mathbb{R}(T_{j}-T_{j-1}) \,dx +
\frac{\kappa}{h_j}\int_\mathbb{R}(T_{j}-T_{j+1}) \,dx +
\frac{\kappa}{h_{j-1}}\int_\mathbb{R}(T_{j-2}-T_{j-1}) \,dx.
\nonumber
\end{eqnarray}
Multiplying \eqref{diff1} by $h_j,$ and summing over $j=1, \dots, N$, we obtain:
\begin{eqnarray}
\begin{split}\label{diff2}
&\sum\limits_{j=1}^N (-1)^{\sigma_j}
h_j\int_\mathbb{R}(T_{j-1}-T_{j})_t
 dx + \sum\limits_{j=1}^N h_j |A_{j-1}-A_j| =\\
&-\sum\limits_{j=1}^N (-1)^{\sigma_j} h_j
\int_\mathbb{R}T_j(1-T_j) dx + \sum\limits_{j=1}^N (-1)^{\sigma_j}
h_{j-1}\int_\mathbb{R}T_{j-1}(1-T_{j-1})dx \frac{h_j}{h_{j-1}}\\
&+ \sum\limits_{j=1}^N (-1)^{\sigma_j} \kappa
\left[\left(1+\frac{h_j}{h_{j-1}}\right)\int_\mathbb{R}(T_{j}-T_{j-1})
dx + \int_\mathbb{R}(T_{j}-T_{j+1}) dx \right.
\\  &+ \left. \int_\mathbb{R}(T_{j-2}-T_{j-1}) dx
\frac{h_j}{h_{j-1}} \right], \\
\end{split}
\end{eqnarray}
where $\sigma_j =0$ if $A_j \geq A_{j-1}$ and $\sigma_j=1$
otherwise. Denote the last sum in \eqref{diff2} by $I_1.$ Using
Lemma~\ref{keylem} and \eqref{keycontrol}, we can estimate:
\begin{eqnarray}
\begin{split}\label{cruin}
& |I_1| \leq C  \kappa \sum\limits_{j=1}^N \left[
\int_\mathbb{R}(T_j-T_{j-1})^2 dx + \int_\mathbb{R}T_j(1-T_j)
\right]\left(1+\underset{j}{max}\left(\frac{h_j}{h_{j-1}}\right)\right)\\
& \leq CH \left( R_2  V^\prime(t) + \kappa V(t) R_1R_2 \right).
\end{split}
\end{eqnarray}
As in the proof of Theorem~\ref{2layer}, we average in $t \in [0,\tau]$ according to \eqref{aver1} and apply
Lemma~\ref{keylem} to estimate the expressions arising from the
time derivatives in \eqref{diff2} in terms of $V$ and $V'$. Taking
into account \eqref{cruin}, after the first integration we obtain,
with $\tau_1 = \frac\tau2-k$ and $\tau_2 = \frac\tau2+k,$
\begin{eqnarray*}
\sum\limits_{j=1}^N (-1)^{\sigma_j} h_j
\left[\int_\mathbb{R}(T_{j-1}(\tau_2)-T_{j}(\tau_2))
 dx-\int_\mathbb{R}(T_{j-1}(\tau_1)-T_{j}(\tau_1))
 dx \right] + \\   (\tau_2-\tau_1)\sum\limits_{j=1}^N h_j |A_{j-1}-A_j| \leq
CHR_2 \left(   V(\tau_2)-V(\tau_1) + (\kappa R_1+1)
\int\limits_{\tau_1}^{\tau_2} V(t)\,dt \right).
\end{eqnarray*}
Applying Lemma~\ref{keylem} once again to the expression on the
left hand side, and integrating twice more, we obtain
\begin{eqnarray*}
\frac{\tau^3}{H}\sum\limits_{j=1}^N h_j |A_{j-1}-A_j| \leq CR_2
\left[ (\kappa R_1+1) \int\limits_0^\tau V(t)G(\frac\tau2, t-
\frac\tau2)\,dt + \right. \\ \left. (\kappa^{-1}+R_1)
\int\limits_0^\tau V(t) h(t,\tau)\,dt + \kappa^{-1}
\int\limits_0^\tau V(t)\,dt \right].
\end{eqnarray*}
Taking into account the bounds on functions $G$ and $h,$ we arrive
at
\begin{equation}\label{fin1}
\frac1\tau\int\limits_0^\tau V(t)\,dt \geq
\frac{C}{HR_2}\sum\limits_{j=1}^N h_j |A_{j-1}-A_j| \left( 1+
\kappa R_1 + \frac{1}{\kappa \tau^2} +\frac{1}{\kappa \tau} +
\frac{R_1}{\tau} \right)^{-1}.
\end{equation}
This completes the proof of Theorem~\ref{main2}. In the situation
where all $h_j \equiv 1,$ and hence $R_1=R_2=2,$ we obtain Theorem~\ref{main1} exactly.
In general, we see that small $h_j$ (which translate to the large coupling constants $\alpha_j$),
or strongly non-uniform size distributions of the neighboring
$h_j,$ weaken the lower bound.


\end{proof}

\section{Numerical Simulations}\label{num}

We carried out numerical simulations for the two layer system
\eqref{twolayer}. Our main interest was in finding the
reaction propagation rate in the established regime, depending on
$A$ and $\kappa.$ Our results are obtained using the
Crank-Nicholson implicit scheme.  The implementation features two
noteworthy modifications. Firstly, the transient wave that occurs
between the initial condition and the eventual steady wavefront
shape was ignored when calculating the front speed.  This was
achieved by using a "$t_{transient}$" parameter, and then ignoring
the wave speed for times $0$ to $t_{transient}$.  Secondly, we
employed a computational trick to expedite simulation.  In a naive
implementation, we would apply the Crank-Nicholson scheme to the
entire field of points that the wave travels through.  Since
the matrices involved in the computation are $n \times n$ where
$n$ is the number of points being considered, considering the
entire field of points that the wave passes through can be
computationally expensive. We avoided this by only applying
Crank-Nicholson to points in the immediate neighborhood of the
front itself.  Then after every iteration we adjusted our
calculation frame to stay centered on the wavefront.

\begin{figure}
\begin{center}
\scalebox{0.75}{\includegraphics{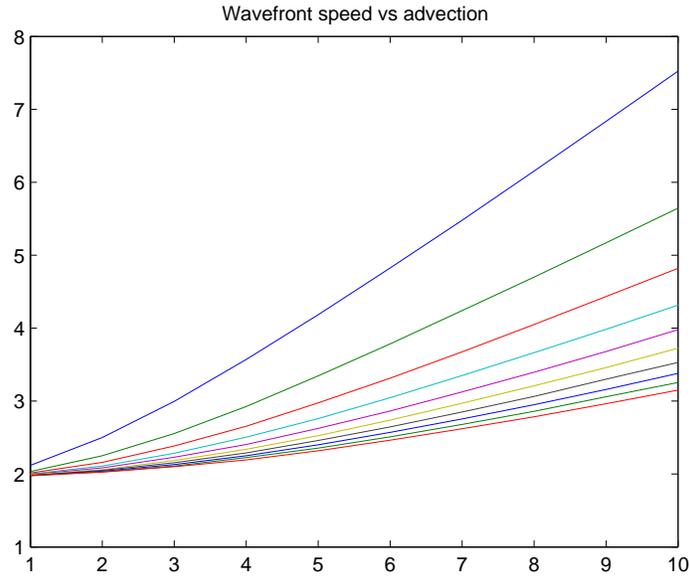}} \caption{The
bulk burning rate versus flow strength $A$ for different values of
the coupling $\kappa$} \label{lin}
\end{center}
\end{figure}

\begin{figure}
\begin{center}
\scalebox{0.75}{\includegraphics{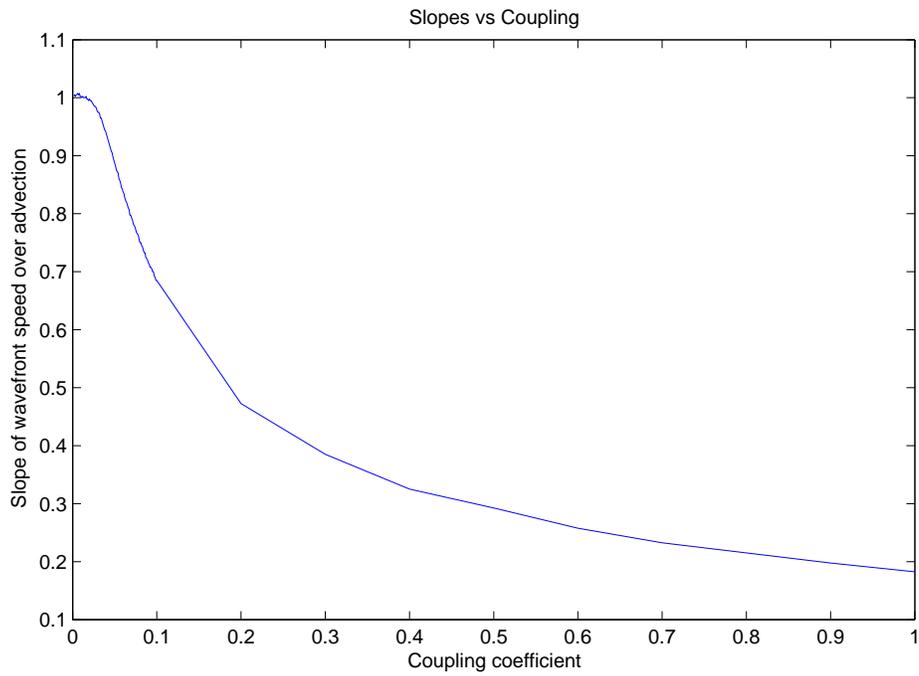}} \caption{The slopes
corresponding to different values of $\kappa$} \label{slopes}
\end{center}
\end{figure}

We found clear linear dependence on the flow strength parameter
$A.$ The Figure~\ref{lin} shows this linear dependence for
different values of the coupling constant $\kappa.$ The constant
$\kappa$ varies from $0.1$ to $1$ with $0.1$ step. The slope
of the graph is monotone decreasing in $\kappa,$ so that the
steeper slopes on Figure~\ref{lin} correspond to smaller values of
$\kappa.$ The intuitive explanation for this effect is that
stronger diffusion between layers mollifies the front stretching
produced by the flow, thus reducing the reaction zone. The graph
on Figure~\ref{slopes} shows the dependence of slope on $\kappa.$
The shape of the graph is in good qualitative agreement with the
bound of Theorem~\ref{2layer}. We note that the results of our
simulations are similar to the results of \cite{VCKRR}, where
combustion in a shear flow was studied in a full PDE setting by
means of a more sophisticated numerical scheme.

{\bf Acknowledgments} This work has been a part of Collaborative
Undergraduate Research Lab (CURL) at the University of Wisconsin
during the 2004-2005 academic year. Support of the NSF VIGRE grant
is gratefully acknowledged.  We thank Carl Edquist, Morgan
Franklin, Jeremy Jacobson, Paul Heideman, Julie Mitchell and John
Vano for interesting discussions and help with this project. We
are especially indebted to Paul Milewski since without his
assistance with the code, the numerical simulations would be still
running. The work of LRM, DA and EB was supported in part by the
NSF VIGRE grant. The work of AK has been partially supported by
the Alfred P. Sloan Research Fellowship and NSF-DMS grant 0314129.

\end{document}